
\magnification=\magstep 1
\baselineskip=15pt     
\parskip=3pt plus1pt minus.5pt
\overfullrule=0pt
\font\hd=cmbx10 scaled\magstep1

\font\hd=cmbx10 scaled\magstep1
\def\num{\global\advance\count10 by 1 \eqno(\the\count10)}

\def\Pfour{{\bf P^4}}
\def\Pthree{{\bf P^3}}

\def\Ptwo{{\bf P^2}}
\def\Pone{{\bf P}^1}
\def\P{{\bf P}}

\def\SS{{\cal S}}
\def\X{{\cal X}}
\def\O{{\cal O}}

\def\T{{\cal T}}
\def\H{{\cal H}}
\def\G{{\cal G}}

\def\F{{\cal F}}
\def\I{{\cal I}}
\def\M{{\cal M}}
\def\N{{\cal N}}
\def\E{{\cal E}}
\def\R{{\cal R}}
\def\K{{\cal K}}
\def\Y{{\cal Y}}

\def\EY{\E(1)|_Y}

\def\PEY{\P(\EY)}

\def\OP{\O_{\PEY}(1)}
\def\OP1{\O_{\Pone}}

\def\Pic{{\rm Pic}}

\def\coker{\mathop{\rm coker}}

\def\im{\mathop{\rm im}}

\def\ext{{\rm Ext}}
\def\tor{{\rm Tor}}
\def\hom{{\rm Hom}}

\def\lhom{{\underline{\hom}}}
\def\lext{{\underline{\ext}}}
\def\ltor{{\underline{\tor}}}
\def\boldz{{\bf Z}}

\def\dual#1{{#1}^{\scriptscriptstyle \vee}}

\def\exact#1#2#3{0\rightarrow#1\rightarrow#2\rightarrow#3\rightarrow0}

\def\mapleft#1{\smash{
  \mathop{\longleftarrow}\limits^{#1}}}
\def\mapright#1{\smash{
  \mathop{\longrightarrow}\limits^{#1}}}
\def\mapsw#1{\smash{
  \mathop{\swarrow}\limits^{#1}}}
\def\mapse#1{\smash{
  \mathop{\searrow}\limits^{#1}}}

\def\mapdown#1{\Big\downarrow
   \rlap{$\vcenter{\hbox{$\scriptstyle#1$}}$}}
\def\mapup#1{\Big\uparrow
   \rlap{$\vcenter{\hbox{$\scriptstyle#1$}}$}}

\centerline{\hd Primitive Calabi-Yau threefolds.}
\medskip
\centerline{\it Mark Gross\footnote{*}{Supported in part by NSF grant
DMS-9400873}}
\medskip
\centerline{November, 1995}
\medskip
\centerline{Department of Mathematics}
\centerline{Cornell University}
\centerline{Ithaca, NY 14853}
\centerline{mgross@math.cornell.edu}
\bigskip
\bigskip
{\hd \S 0. Introduction.}

A Calabi-Yau threefold is a complex projective threefold $X$ (possibly with
some
suitable class of singularities, say terminal or canonical) with $\omega_X\cong
\O_X$ and $h^1(\O_X)=h^2(\O_X)=0$. One of the fundamental gaps in the
classification of algebraic threefolds is the lack of understanding of
Calabi-Yau threefolds. Here I will try to set forth a program to bring the
morass of thousands of examples of Calabi-Yaus under control.

The ideas here go back to the papers of Friedman [2] and Reid [25]. Friedman
studied smoothability of Calabi-Yau threefolds with ordinary double points.
Based on these results, Reid conjectured that there could perhaps be a single
irreducible moduli space of (non-K\"ahler) Calabi-Yau threefolds, such that any
Calabi-Yau threefold is the small resolution of a degeneration of this family
to something with ordinary double points. So one can think of all the chaos of
the algebraic examples as simply being ``boundary phenomena'' for the moduli
space of this ur-Calabi-Yau. I suspect the most difficult part of this
conjecture, often known as Reid's fantasy, will be passing from algebraic to
non-algebraic threefolds. We do not understand how to deal with non-K\"ahler
Calabi-Yau threefolds or find non-algebraic contractions.

Unlike in the K3 case, where it is possible to deform an algebraic K3 surface
to a non-algebraic one, the deformation of a projective Calabi-Yau threefold,
even singular, is still projective. So it makes sense to insist on staying
within the projective category. Reid's picture given above needs to be modified
if we restrict attention to projective threefolds; indeed, an algebraic
Calabi-Yau of Picard number 1 has no algebraic birational contractions, and so
cannot be the resolution of another algebraic Calabi-Yau threefold. Since there
are many examples of Calabi-Yau threefolds with Picard number 1, we cannot hope
that all Calabi-Yau threefolds arise as resolutions of degenerations
of a single family of algebraic threefolds. Nevertheless, a somewhat weaker
picture seems reasonable.

As proposed in [1], [3] and [4], we can think of the moduli of Calabi-Yaus as
forming a giant web, a directed graph where each node is a deformation class of
Calabi-Yau threefolds. We draw an arrow $\M_1\rightarrow \M_2$ if for the
general element
$\tilde X$ of deformation class $\M_1$, there is a birational contraction
morphism $\pi:\tilde X\rightarrow X$ and a flat family $\X\rightarrow
(\Delta,0)$ such that $\X_0\cong X$ and $\X_t\in \M_2$ for general
$t\in\Delta$.
For example, let $\M_Q$ be the moduli space of smooth quintics in $\P^4$,
$\M_D$ the moduli space of double covers of $\P^3$ branched over smooth octics.
Let $T$ be the blowup of a quintic with a triple point; let $\M_T$
be the moduli space of such $T$. We have two contraction morphisms
$\pi_1:T\rightarrow T_1$, $\pi_2:T\rightarrow T_2$ with $\pi_1$ the contraction
of the exceptional cubic surface to $T_1$, a quintic with a triple point,
and $\pi_2$ the Stein factorization of the projection $T\rightarrow\Pthree$
from the triple point. $T_1$ can be smoothed to a smooth quintic, and $T_2$
can be smoothed to a double cover of $\Pthree$ branched over a smooth octic.
Thus a tiny portion of our web will be
$$\matrix{&&\M_T&&\cr
&\swarrow&&\searrow&\cr
\M_D&&&&\M_Q\cr}$$

By taking deformation classes of all simply-connected Calabi-Yau threefolds, we
get an enormous (perhaps infinite) graph. So one question that
immediately comes up is often thought of as one version of Reid's fantasy:

\proclaim The Connectedness Conjecture 0.1. The graph of simply connected
Calabi-Yaus is connected.

The evidence for this is strictly experimental at this point: large classes of
examples have been connected up, ``by hand,'' e.g. [4], where moduli of
Calabi-Yau complete intersections in products of projective spaces are
connected
up. Work in progress with T.-M. Chiang, B. Greene and Y. Kanter has connected
up
all known examples of Calabi-Yaus in weighted $\Pfour$.

Clearly, the nodes at the bottom of this Calabi-Yau graph, such as $\M_Q$ and
$M_D$, will be the Calabi-Yaus which play the role of Reid's ur-Calabi-Yau in
the
algebraic setting. I call a Calabi-Yau whose deformation class is a node of the
graph with no outgoing arrows a  {\it primitive} Calabi-Yau. More formally,

\proclaim Definition 0.2. A non-singular Calabi-Yau threefold $\tilde X$ is
{\it primitive} if there is no birational contraction $\tilde
X\rightarrow X$ with $X$ smoothable to a Calabi-Yau threefold which is not
deformation equivalent to $\tilde X$.

So one approach to proving the connectedness conjecture would be to classify
primitive Calabi-Yaus and then to find some paths through the Calabi-Yau graph
connecting the primitive Calabi-Yaus together.

Beyond the connectedness conjecture, we might learn a lot more about the number
of families of Calabi-Yaus as a whole if we can understand the class of
primitive
Calabi-Yaus. Four further questions of varying strength are

\proclaim Questions 0.3. \item{(1)} Are there are a finite number of flat
families
$\X_i\rightarrow\SS_i$ of non-singular Calabi-Yau threefolds such that for any
Calabi-Yau $\tilde X$, there is a flat family $\X\rightarrow\Delta$ with $\X_t$
in one of the families $\X_i\rightarrow\SS_i$ for $t\not=0$ and $\X_0$
birational to $\tilde X$?
\item{(2)} The same question, except we insist that $\tilde X$ be a crepant
resolution of $\X_0$.
\item{(3)} Are there a finite number of flat families $\X_i\rightarrow\SS_i$
of Calabi-Yau threefolds with canonical singularities such that any
Calabi-Yau $\tilde X$ is birational to some member $X$ of one of these
families?
\item{(4)} Same question as (3), except we insist that $\tilde X$ be a crepant
resolution of $X$.

Intuitively, these should follow if there are only a finite number of families
of primitive Calabi-Yau threefolds. Unfortunately, none of them do immediately,
though (1) follows if fourfold flops and ${\bf Q}$-factorializations exist. We
will discuss these questions in \S 3, and in particular the current technical
obstacles involved in answering these questions, even under the assumption
that there are only a finite number of families of primiticve Calabi-Yaus. Note
that (3) implies there are only a finite number of families up to birational
equivalence, which has been proven in the case of elliptic Calabi-Yau
threefolds
in [6]. (4) implies the stronger result that there are only a finite number of
families up to biregular equivalence, since a Calabi-Yau with canonical
singularities has only a finite number of crepant resolutions.

So the main philosophy espoused by this paper is that we should attempt to
classify primitive Calabi-Yaus, and that this should be easier than classifying
all Calabi-Yaus. Is there any hope for such a classification?
Well, obviously any Picard number 1 Calabi-Yau is
primitive, and we do not have a classification of Picard number 1 Calabi-Yaus,
but
at least I only know of about 30 such threefolds, as opposed to an order of
10000
currently known examples of Calabi-Yaus in general. Suppose we could understand
Picard number one Calabi-Yaus. The hope then would be that there are very few
primitive Calabi-Yaus with Picard number greater than one. There certainly are
some: the bidegree $(3,3)$ hypersurface in
$\Ptwo\times\Ptwo$ furnishes an example, since it has no algebraic
contractions, and there are
other, less trivial examples. However, we shall give some evidence in \S 2 that
we
should not expect many such examples.

How do we attack this classification problem? Given a Calabi-Yau $\tilde X$,
we presumably need to understand birational contractions of $\tilde X$,
$\pi:\tilde X\rightarrow X$, and smoothability of $X$. Now $X$ will have
canonical singularities, and these can be very complicated. So we can't really
hope to completely answer the question of when $X$ is smoothable. However, we
can answer this question if we assume $\pi$ is a primitive contraction, i.e.
$\pi$
cannot be factored in the algebraic category. We have already begun this study
in [8], where type I and II contractions were studied (small contractions
and divisorial contractions to points respectively). However, type III
contractions
still need to be considered. These are contractions which contract a divisor to
a
curve. Technically, this is a harder case to deal with, and we devote \S 1
to it. The final result, however, is remarkably similar to Theorem 5.8 of [8].
We have

\proclaim Theorem 0.4. Let $\pi:\tilde X\rightarrow X$ be a primitive type III
contraction, contracting a divisor $E$ to a curve $C$. Then $X$ is smoothable
unless $C\cong\Pone$ and $E^3=7$ or $8$.

This follows from Theorems 1.3 and 1.7.

Recall that on a Calabi-Yau threefold, $E^3$ is the self-intersection of the
canonical class of $E$. Thus, if $E$ is a normal surface, $C\cong\Pone$, then
$E^3=8$ is the case that $E$ is a minimal scroll and $E^3=7$ is the case that
$E$ is a minimal scroll blown up in one point.

Combining this with the results of [8], we see that if $\tilde X$ is a
primitive Calabi-Yau and $\pi:\tilde X\rightarrow X$ is a primitive
contraction, then $\pi$ is either a contraction of a single $\Pone$ with normal
bundle $\O_{\Pone}(-1)\oplus\O_{\Pone}(-1)$, or a contraction of $\Ptwo$, a
minimal ruled surface over $\Pone$, such a surface blown up in one point, or
a non-normal surface with $E^3=7$. This in fact turns out to be a very strong
restriction, principally since then the exceptional divisor $E$ of any
primitive divisorial contraction on a primitive Calabi-Yau satisfies $c_2.E<0$.
We shall take this up in
\S 2. There we will derive some immediate combinatorial consequences of
these results, and speculate what still needs to be done. I believe that a more
detailed combinatorial analysis will yield far stronger results then are
obtained
in \S 2. I hope to treat this approach in a future paper.

 Finally, we note that connecting together moduli spaces of Calabi-Yau
threefolds appears to be playing an important role in physics. In particular,
in [5], a physical explanation has been given for transitions between moduli
spaces. So far only degenerations involving ordinary double points have been
studied, but there is apparently no physical reason for restricting attention
to
such cases.

I would like to thank P.M.H. Wilson for many helpful conversations during the
preparation of this paper.

{\hd \S 1. Type III Contractions.}

\proclaim Theorem 1.1. Let $\pi:\tilde X\rightarrow X$ be a primitive type III
contraction of a non-singular Calabi-Yau threefold $\tilde X$, contracting an
exceptional divisor $E$ to a curve $C$. Then
\item{(a)} $C$ is a non-singular curve.
\item{(b)} $\pi:E\rightarrow C$ is a conic bundle over $C$, and each fibre is
either a non-singular conic, a union of two lines meeting at a point, or a
doubled line. If the general fibre is a non-singular curve, then $E$ is normal.
In this case, the singularities which appear on $E$ are $A_n$ ($n\ge 0$)
singularities at the singular point of a reducible reduced fibre, or two $A_1$
singularities on a non-reduced fibre.

Proof. For (a), see [34]; for (b), see [30], Theorem 2.2, keeping in mind
[31]. $\bullet$

\proclaim Proposition 1.2. Let $\pi:\tilde X\rightarrow X$ be a primitive type
III contraction of a non-singular Calabi-Yau threefold $\tilde X$, contracting
a divisor $E$ to a curve $C$. Let $\tilde E$ be the normalization of $E$,
$f:\tilde E\rightarrow\tilde X$ the induced map, and $\tilde E\rightarrow\tilde
C\rightarrow C$ the Stein factorization. Then the image of the natural map
$Def(f)\rightarrow Def(\tilde X)$ has codimension $\ge p_a(\tilde C)$.

Proof. If $E$ is already normal, then this is the result of [17],
Prop. 6.5. If $E$ is not normal, then the fibres of $\pi:\tilde
E\rightarrow C$ are line pairs or doubled lines, and $\tilde\pi:\tilde
E \rightarrow \tilde C$ is $\Pone$ fibration with a section. Thus $\tilde E$
is a non-singular scroll. We proceed as in [31].

Define $N_f$ by the exact sequence
$$\exact{\T_{\tilde E}}{f^*\T_{\tilde X}}{N_f}.$$
$N_f$ is torsion free, and fails to be locally free precisely at the inverse
images of the pinch points of $E$. Thus $\dual{\dual{N_f}}$ is locally free and
$c_1(N_f)=c_1(\dual{\dual{N_f}})$, showing that $\dual{\dual{N_f}}\cong
\omega_{\tilde E}$.

If $T^1_f$ is the tangent space to $Def(f)$, we have an exact sequence by [23]
$$H^0(\T_{\tilde X})\oplus H^0(\T_{\tilde E}) \rightarrow H^0(f^*\T_{\tilde X})
\rightarrow T^1_f \rightarrow H^1(\T_{\tilde X})\oplus H^1(\T_{\tilde E})
\rightarrow H^1(f^*\T_{\tilde X}).$$ This induces an exact sequence
$$T^1_f\mapright{\alpha} H^1(\T_{\tilde X}) \mapright{\beta} H^1(N_f)$$
where $\alpha$ is the differential of the map $Def(f)\rightarrow Def(\tilde X)$
and $\beta$ is induced by $H^1(\T_{\tilde X})\rightarrow H^1(f_*f^*\T_{\tilde
X})= H^1(f^*\T_{\tilde X})\rightarrow H^1(N_f)$. Since $H^0(N_f)\subseteq
H^0(\omega_{\tilde E})=0$, $\alpha$ is injective. Thus, if $\dim\im\beta\ge
p_a(\tilde C)$, $codim(\im(Def(f)\rightarrow Def(\tilde X)))\ge p_a(\tilde C)$.
The composed map
$$H^1(\T_{\tilde X})\rightarrow H^1(N_f)\rightarrow H^1(\dual{\dual{N_f}})
=H^1(\omega_{\tilde E})$$ is the natural map $H^1(\Omega^2_{\tilde X})
\rightarrow H^1(\Omega^2_{\tilde E})$. Since $h^1(\omega_{\tilde E})
=p_a(\tilde C)$, it is sufficient to show that $H^1(\Omega^2_{\tilde X})
\rightarrow H^1(\Omega^2_{\tilde E})$ is surjective, or equivalently by Hodge
theory, that
$H^2(\Omega^1_{\tilde X})\rightarrow H^2(\Omega^1_{\tilde E})$ is surjective.

To show surjectivity of this map, note that we have exact sequences
$$0\rightarrow \F_1\rightarrow \Omega^1_{\tilde X}\rightarrow
f_*f^*\Omega^1_{\tilde X}\rightarrow \F_2\rightarrow 0\leqno{(1.2.1)}$$
where $\F_2$ has support on the singular curve of $E$, and
$$0\rightarrow \F_3\rightarrow f^*\Omega^1_{\tilde X}\rightarrow
\Omega^1_{\tilde E}\rightarrow \F_4\rightarrow 0\leqno{(1.2.2)}$$
where $\F_4$ has support on the pinch points of $\tilde E$. (1.2.2) shows that
$H^2(f^*\Omega^1_{\tilde X})\rightarrow H^2(\Omega^1_{\tilde E})$
is surjective, and (1.2.1) shows that $H^2(\Omega^1_{\tilde X})\rightarrow
H^2(f^*\Omega^1_{\tilde X})$ is surjective if $H^3(\F_1)=0$. Now $H^3(\F_1)
\cong \dual{H^0(\dual{\F_1})}$ by Serre duality, and there is an injection
$\pi_*\dual{\F_1}\rightarrow \T_X$ since $\pi_*\dual{\F_1}$ is torsion-free,
$\pi_*\dual{\F_1}$  and $\T_X$ coincide off a codimension 2 set, and $\T_X$ is
reflexive. Thus $H^0(\dual{\F_1})\subseteq H^0(\T_X)=0$ by [11], 8.6. We
conclude
that $H^2(\Omega^1_{\tilde X})\rightarrow H^2(\Omega^1_{\tilde E})$ is
surjective, and the proposition follows.
$\bullet$

\proclaim Theorem 1.3. Let $\pi:\tilde X\rightarrow X$ be a primitive type
III contraction of a non-singular Calabi-Yau threefold $\tilde X$, contracting
a divisor $E$ to a curve $C$. If $p_a(C)\ge 1$, then $X$ is smoothable.

Proof: Let $f:\tilde\X\rightarrow\Delta$ be a deformation
of $\tilde X$ over a contractible base $\Delta$ which is sufficiently general,
so
that the exceptional divisor
$E$ does not deform to general $\tilde\X_t$, $t\in\Delta$. By Proposition 1.2,
such exists. The contraction $\pi:\tilde X\rightarrow X$ yields a
contraction $\tilde\X\rightarrow \X$ over $\Delta$ extending $\pi$. The general
contraction $\tilde \X_t\rightarrow\X_t$ is then a small primitive contraction,
and by [8], Proposition 5.1, $\X_t$ is smoothable unless $\X_t$ has exactly one
ordinary double point. Thus it is enough to show that for general
$t$, the singular locus of $\X_t$ is not exactly one ordinary double point.
Note
we are not excluding the possibility that $\tilde \X_t=\X_t$ for $t$ general,
but
in this case we are done. (This can happen if $p_a(C)=1$; see [30].)

First suppose that $E$ is normal and $\pi:E\rightarrow C$ has a singular fibre.
By Theorem 1.1 b), this singular fibre is either two $\Pone$'s or a doubled
$\Pone$. Let $Z$ be the homology class of either one of these $\Pone$'s in
the first case or of the reduced fibre in the second case. If, for
general $t\in\Delta$, $\X_t$ had only one ODP, then $\tilde\X_t$ would contain
precisely one
$\Pone$ with normal bundle $\O_{\Pone}(-1)\oplus\O_{\Pone}(-1)$ in the
homology class of $Z$. This then implies that if we deform the complex
structure of $\tilde X$ to a generic almost complex structure, the homology
class $Z$ represents only one pseudo-holomorphic curve. But this contradicts
[34] Lemma 4.1. Thus for general $t$, $\X_t$ has worse singularities than one
ODP.

Now suppose
$E$ is either normal and
$\pi:E\rightarrow C$ has no singular fibres, or $E$ is non-normal. If, in the
notation of Proposition 1.2,
$p_a(\tilde C)=1$,
then by [30] Proposition 4.4 and [31], $\tilde\X_t\rightarrow\X_t$ is an
isomorphism for general $t\in\Delta$; thus $X$ is smoothable. If $p_a(\tilde
C)\ge 2$, let $l$ be the class of a fibre of $\tilde E\rightarrow \tilde C$,
and
let
$\tilde\X\rightarrow Def(\tilde X)$ be the Kuranishi family of $\tilde X$. Let
$S\rightarrow Def(\tilde X)$ be the irreducible component of the relative
Douady
space of $\tilde\X\rightarrow Def(\tilde X)$ corresponding to deformations of
$l$. Now
$\dim S\ge \chi(N_{l/\tilde \X})=\dim Def(\tilde X)$, but on the other hand, by
Proposition 1.2, the locus in $Def(\tilde X)$ where
$S\rightarrow Def(\tilde X)$ has one dimensional fibres is codimension at least
2. Thus $S\rightarrow Def(\tilde X)$ must be surjective. Suppose that this map
is generically one-to-one. Then $S_{red}\rightarrow Def(\tilde X)$ is a
birational map, and so the one dimensional fibres must be rational. This
contradicts $p_a(\tilde C)\ge 2$. Thus again, for general $t\in Def(\tilde X)$,
$\tilde\X_t\rightarrow \X_t$ is a small contraction contracting at least two
curves.
$\bullet$

We now focus on the necessary local deformation-theoretic calculations for the
case that $C\cong\Pone$; in this case, the exceptional divisor $E$ always
deforms with $\tilde X$.

\proclaim Theorem 1.4. Let $\pi:\tilde X\rightarrow X$ be a primitive type
III contraction of a non-singular Calabi-Yau threefold $\tilde X$, contracting
a divisor $E$ to a curve $C\cong \Pone$. Suppose furthermore that $\tilde X$
is general in its moduli. Then
\item{(a)} If $E^3\not= 7$, then $E$ is a normal surface.
\item{(b)} $\I_C/\I_C^2$ is a locally free sheaf on $C$ of rank 3 and degree
$E^3-2$.
\item{(c)} Let ${\bf T}^1=\lext^1_{\O_X}(\Omega^1_X,\O_X)$.
If $E$ is normal, then  there
are exact sequences
$$\exact{R^1\pi_*\lhom_{\O_{\tilde X}}(\pi^*\Omega^1_X,\O_{\tilde X})}
{{\bf T}^1}{\pi_*\lext^1_{\O_{\tilde X}}(\pi^*\Omega^1_X,\O_{\tilde
X})}\leqno{(1.4.1)}$$ and
$$\exact{\F}{\pi_*\lext^1_{\O_{\tilde X}}(\pi^*\Omega^1_X,\O_{\tilde X})}{
\O_{C}(4-E^3)}\leqno{(1.4.2)}$$ where both
$R^1\pi_*\lhom_{\O_{\tilde X}}(\pi^*\Omega^1_X,\O_{\tilde X})$ and
$\F$ are
sheaves of finite length, with
$$length(\F)\ge 8-E^3.$$

Proof. (a) If $E$ is not normal, then following the notation of Proposition
1.2, we can assume that $p_a(\tilde C)=0$ since $\tilde X$ is general in its
moduli. Thus $\tilde C\rightarrow C$ is a double cover branched in precisely
two points. By [34], Prop. 3.2, the map $\tilde
E\rightarrow X$ has precisely two pinch points. As
pointed out at the end of the proof of [34], Prop. 4.2, this implies that
$E^3=7$.  The same calculation as in [8], Theorem 5.2 also shows that $E^3=7$.

(b) Let $\I_C$ be the ideal sheaf of $C\subseteq X$. We first need to
understand $\I_C/\I_C^2$. In a neighborhood of a point of $C$, we can embed
$X\subseteq Y$, with $Y$ smooth of dimension 4, since $X$ has only hypersurface
singularities. Since $X$ is singular along $C$, $\I_{X/Y}\subseteq \I_{C/Y}^2$,
and so $\I_C/\I_C^2= \I_{C/Y}/(\I_{C/Y}^2+\I_{X/Y})= \I_{C/Y}/
\I_{C/Y}^2$. Thus $\I_C/\I_C^2$ is locally free of rank 3.
Furthermore, since $\pi$ is a primitive contraction, it must be the blow-up
of $X$ along $C$, and so locally $\tilde X\subseteq\tilde Y$ where $\tilde Y$
is the blow-up of $Y$ along $C$. Let $F\subseteq \tilde Y$ be the exceptional
locus of the blow-up of $Y$ along $C$; $F$ is isomorphic to the $\Ptwo$-bundle
over $C$ given by $\P(\I_C/\I_C^2)$, and $\I_{F/\tilde Y}/\I_{F/\tilde Y}^2
=\O_F(1)$. Thus
$$\eqalign{\I_C/\I_C^2&=\pi_*\O_F(1)\cr
&=\pi_*\O_F(1)|_E\cr
&=\pi_* \I_{E/\tilde X}/\I_{E/\tilde X}^2\cr
&=\pi_*\omega_E^{-1}.\cr}$$
It
is easy to see that $R^1\pi_*\omega_E^{-1}=0$, so $\chi(\omega_E^{-1})
=\chi(\I_C/\I_C^2)=\deg(\I_C/\I_C^2)+3$, as $C\cong \Pone$. To compute
$\chi(\omega_E^{-1})$, recall that $E$ has du Val singularities from Theorem
1.1 (b). If
$h:\tilde E\rightarrow E$ is a minimal resolution, then $h^*\omega_E^{-1}
=\omega_{\tilde E}^{-1}$ and so $\chi(\omega_E^{-1})=\chi(\omega_{\tilde
E}^{-1})-h^0(R^1h_*\omega_{\tilde E}^{-1})$. Furthermore $R^1h_*\omega_{\tilde
E}^{-1}=0$, so $\chi(\omega_E^{-1})=\chi(\omega_{\tilde E}^{-1})=K^2_{\tilde
E}+1$, by Riemann-Roch, where $K_{\tilde E}$ is the canonical class on $\tilde
E$. But $K_{\tilde E}^2=K^2_E=E^3$. We conclude that
$$\deg\I_C/\I_C^2=E^3-2.$$

(c) By [9], II 5.10, there
is an isomorphism
$${\bf R}\pi_*{\bf R}\lhom_{\O_{\tilde X}}(\pi^*\Omega^1_X,\O_{\tilde X})
\cong {\bf R}\lhom_{\O_X}(\Omega^1_X,{\bf R}\pi_*\O_{\tilde X})\cong {\bf
R}\lhom_{\O_X}(\Omega^1_X,\O_X).\leqno{(1.4.3)}$$
The latter isomorphism follows from the fact that $X$ has rational
singularities. This yields a spectral sequence
$$R^p\pi_*\lext^q_{\O_{\tilde X}}(\pi^*\Omega^1_X,\O_{\tilde X})
\Rightarrow \lext_{\O_X}^n(\Omega^1_X,\O_X)$$
which gives the exact sequence
$$\exact{R^1\pi_*\lhom_{\O_{\tilde X}}(\pi^*\Omega^1_X,\O_{\tilde X})}
{\lext^1_{\O_X}(\Omega^1_X,\O_X)}{\pi_*\lext^1_{\O_{\tilde
X}}(\pi^*\Omega^1_X,\O_{\tilde X})}.$$ This yields (1.4.1).
To calculate $\pi_*\lext^1_{\O_{\tilde X}}(\pi^*\Omega^1_X,\O_{\tilde X})$, we
use the exact sequence
$$0\mapright{} \G \mapright{} \pi^*\Omega^1_X \mapright{} \Omega^1_{\tilde X}
\mapright{\phi} \Omega^1_{\tilde X/X} \mapright{} 0\leqno{(1.4.4)}$$
where $\G$ is a sheaf supported on $E$. This yields an exact sequence
$$\eqalign{
0=\lhom_{\O_{\tilde X}}(\G,\O_{\tilde X})&\rightarrow
\lext^1_{\O_{\tilde X}}(\ker\phi,\O_{\tilde X})\rightarrow \lext^1_{\O_{\tilde
X}}(\pi^*\Omega^1_X,\O_{\tilde X})\cr
&\rightarrow \lext^1_{\O_{\tilde
X}}(\G,\O_{\tilde X})\rightarrow \lext^2_{\O_{\tilde
X}}(\ker\phi,\O_{\tilde X})\cr}\leqno{(1.4.5)}$$
and since $\lext^i_{\O_{\tilde X}}(\Omega^1_{\tilde X},\O_{\tilde X})=0$ for
$i>0$, we obtain isomorphisms
$$\lext^i_{\O_{\tilde X}}(\ker\phi,\O_{\tilde X})\cong\lext^{i+1}_{\O_{\tilde
X}}(\Omega^1_{\tilde X/X},\O_{\tilde X})$$
for $i>0$.

To compute $\Omega^1_{\tilde X/X}$, first note that $\Omega^1_{\tilde X/X}
\otimes_{\O_{\tilde X}}\O_E\cong \Omega^1_{E/C}$. Thus, if $\Omega^1_{\tilde
X/X}$ is already an $\O_E$-module, then $\Omega^1_{\tilde X/X}\cong
\Omega^1_{E/C}$. To see that $\Omega^1_{\tilde X/X}$ is an $\O_E$-module, again
locally embed $X$ in a smooth hypersurface $Y$, and let $\tilde Y\rightarrow Y$
be the blow-up of $Y$ along $C$, with exceptional divisor $F$. Then
$\Omega^1_{\tilde Y/Y}$ is an $\O_F$-module, as can be explicitly calculated,
and
thus
$\Omega^1_{\tilde X/Y} =\Omega^1_{\tilde X/X}$ is an $\O_E$-module.

Now the change of rings spectral sequence ([26], Theorem 11.66) gives
$$\lext^i_{\O_E}(\Omega^1_{E/C},\lext^j_{\O_{\tilde X}}(\O_E,\O_{\tilde X}))
\Rightarrow \lext^n_{\O_{\tilde X}}(\Omega^1_{E/C},\O_{\tilde X}),$$
and from
$$\exact{\O_{\tilde X}(-E)}{\O_{\tilde X}}{\O_E}$$
we obtain
$$\lext^j_{\O_{\tilde X}}(\O_E,\O_{\tilde X})\cong \cases{
\O_E(E)\cong \omega_E&if $j=1$;\cr
0&if $j\not=1$\cr}$$
so
$\lext^{i+1}_{\O_{\tilde
X}}(\Omega^1_{\tilde X/X},\O_{\tilde X})\cong \lext^i_{\O_E}(\Omega^1_{E/C},
\omega_E)$. From
$$\exact{\pi^*\Omega^1_C}{\Omega^1_E}{\Omega^1_{E/C}},$$
we obtain
$$\lhom_{\O_E}(\Omega^1_E,\omega_E)\mapright{\pi_*}
\pi^*\T_C\otimes\omega_E\mapright{}\lext^1_{\O_E}(\Omega^1_{E/C},\omega_E)
\mapright{} \lext^1_{\O_E}(\Omega^1_E,\omega_E)\mapright{} 0, \leqno{(1.4.6)}$$
and
$\lext^2_{\O_E}(\Omega^1_{E/C},\omega_E)\cong\lext^2_{\O_E}(\Omega^1_E,\omega_E)
\cong 0$ since $E$ has only hypersurface singularities. Thus (1.4.5) yields
$$\exact{\lext^1_{\O_E}(\Omega^1_{E/C},\omega_E)}{\lext^1_{\O_{\tilde X}}(
\pi^*\Omega^1_X,\O_{\tilde X})}{\lext^1_{\O_{\tilde X}}(\G,\O_{\tilde X})}.
\leqno{(1.4.7)}$$
In (1.4.6)
$\pi_*$ fails to be surjective precisely where $\pi:E\rightarrow C$ is
not smooth, and thus $\coker\pi_*$ has support precisely on this locus. We now
use
the description of the singularities of $E$ from Theorem 1.1 b).
It is well-known
that at a point
$P$ of
$E$ which is an
$A_n$ singularity of
$E$,
$$length_P(\lext^1_{\O_E}(\Omega^1_E,\omega_E))=n.$$
Thus, if $\pi:E\rightarrow C$ has a fibre consisting of two $\Pone$'s meeting
at
a point $P\in E$ with $P$ an $A_n$ singularity ($n\ge 0$) on $E$, then
$$length_P(\lext^1_{\O_E}(\Omega^1_{E/C},\omega_E))\ge n+1.$$
If, on the other hand, $\pi:E\rightarrow C$ has a fibre which is a double line
$l$, then  $E$ is locally of the form, in $\Ptwo\times{\bf A}^1$,
$$x^2+tyz=0,$$
with $\pi$ given by $(x,y,z,t)\mapsto t$. A calculation shows that
$\coker\pi_*$ is an invertible sheaf on the reduced fibre, $\Pone$, and thus
must be $\O_{\Pone}\otimes\omega_E\cong \O_{\Pone}(-1)$. This gives an exact
sequence
$$\exact{\O_{\Pone}(-1)}{\lext^1_{\O_E}(\Omega^1_{E/C},\omega_E)}
{\lext^1_{\O_E}(\Omega^1_E,\omega_E)}.$$
Since $E$ has two $A_1$ singularities on this fibre,
$\pi_*\lext^1_{\O_E}(\Omega^1_{E/C},\omega_E)$ has length 2 at this point.
Adding up all contributions, we now see that
$$length(\pi_*\lext^1_{\O_E}(\Omega^1_{E/C},\omega_E))\ge rank\Pic(\tilde
E)-2,$$
where $\tilde E$ is a minimal resolution of $E$. This in turn is
$8-K^2_{\tilde E}=8-K_E^2=8-E^3$.
Also, we see that $R^1\pi_*\lext^1_{\O_E}(\Omega^1_{E/C},\omega_E)=0$. Thus,
applying $\pi_*$ to (1.4.7), we obtain
$$\exact{\F}{\pi_*\lext^1_{\O_{\tilde X}}(\pi^*\Omega^1_X,\O_{\tilde X})}
{\pi_*\lext^1_{\O_{\tilde X}}(\G,\O_{\tilde X})},$$
where $\F=\pi_*\lext^1_{\O_E}(\Omega^1_{E/C},\omega_E)$ is a sheaf of finite
length with length at least $8-E^3$.
Thus to demonstrate (1.4.2), we only need to calculate $\G$. Again,
locally embed
$X\subseteq Y$,
$Y$ smooth and dimension 4. We have an exact sequence
$$0\mapright{}\pi^*\O_Y(-X)|_X\mapright{\psi}\pi^*\Omega^1_Y|_X
\mapright{}\pi^*\Omega^1_X \mapright{} 0.$$
A local calculation shows that $\psi$ defines a section of
$\pi^*(\Omega^1_Y|_X\otimes \O_Y(X))$ which vanishes with order one along $E$.
Thus the torsion part of $\coker \psi$ is an invertible $\O_E$-module.
Since the torsion part of $\pi^*\Omega^1_X$ is precisely $\G$, this shows that
$\G$ is an invertible $\O_E$-module.

To determine which invertible $\O_E$-module, we tensor (1.4.4) with $\O_E$.
This yields exact sequences
$$\eqalign{\ltor_1^{\O_{\tilde X}}(\ker\phi,\O_E)
\rightarrow \G  \rightarrow \pi^*\Omega^1_X|_E \rightarrow (\ker\phi)|_E
\rightarrow & 0,\cr
0\rightarrow \ltor_1^{\O_{\tilde X}}(\Omega^1_{E/C},\O_E)\rightarrow
(\ker\phi)|_E \rightarrow \Omega^1_{\tilde X}|_E\rightarrow \Omega^1_{E/C}
\rightarrow & 0\cr}\leqno{(1.4.8)}$$
and an isomorphism
$$\ltor_1^{\O_{\tilde X}}(\ker\phi,\O_E)\cong\ltor_2^{\O_{\tilde
X}}(\Omega^1_{E/C},\O_E)=0,\leqno{(1.4.9)}$$
the latter equality from
$$\exact{\O_{\tilde X}(-E)}{\O_{\tilde X}}{\O_E}.$$
This exact sequence also yields
$$0\mapright{} \ltor_1^{\O_{\tilde X}}(\Omega^1_{E/C},\O_E)
\mapright{\cong} \Omega^1_{E/C}\otimes\omega_E^{-1}
\mapright{} \Omega^1_{E/C} \mapright{\cong} \Omega^1_{E/C} \mapright{} 0.
\leqno{(1.4.10)}$$
We also have an exact sequence
$$\exact{\I_C/\I_C^2}{\Omega^1_X|_C}{\Omega^1_C},$$
injectivity on the left again following from a local calculation showing
generic injectivity. Pulling this back to $E$ gives
$$\exact{\pi^*(\I_C/\I_C^2)}{(\pi^*\Omega^1_X)|_E}{\pi^*\Omega^1_C}.$$
Thus
$$c_1((\pi^*\Omega^1_X)|_E)=\pi^*\O_C(E^3-4),\leqno{(1.4.11)}$$ from (b).
(Here, by $c_1$, we really mean $c_1$ on the non-singular part of $E$, but we
abuse notation and write sheaves on $E$.) From the second sequence of (1.4.8)
and
(1.4.10),
$$\eqalign{c_1((\ker\phi)|_E) &= c_1(\Omega^1_{E/C}\otimes \omega_E^{-1})
\otimes c_1(\Omega^1_{\tilde X}|_E)\otimes c_1(\Omega^1_{E/C})^{-1}\cr
&=\omega_E^{-1}.\cr}$$
Finally, from the first sequence of (1.4.8) and (1.4.11), we see that
$c_1(\G)=\pi^*\O_C(E^3-4)\otimes\omega_E$. Since $\G$ is invertible,
$$\G=\pi^*\O_C(E^3-4)\otimes\omega_E.$$
Thus
$$\eqalign{\pi_*\lext^1_{\O_{\tilde X}}(\G,\O_{\tilde X})
&\cong \pi_*\lhom_{\O_E}(\G,\omega_E)\cr
&\cong \pi_*\pi^*\O_C(4-E^3)\cr
&\cong\O_C(4-E^3).\cr}$$
This yields (1.4.2).

We have now shown everything in part (c) except that $R^1\pi_*\lhom(\pi^*
\Omega^1_X,\O_{\tilde X})$ is of finite length. Now the singularities of $X$
along $C$ are generically $A_1$, and it is well known that at those points of
$C$ which are $A_1$ singularities,
${\bf T}^1$ is an invertible $\O_C$-module. Thus from (1.4.1) and (1.4.2),
$R^1\pi_*\lhom(\pi^*\Omega^1_X,\O_{\tilde X})$ must be supported on dissident
points.
$\bullet$

{\it Example 1.5.} The above proposition, as well as Theorem 5.2 of [8], leaves
open the possibility of a non-normal exceptional divisor $E$, $E$ rational,
with $E^3=7$. We give an example of such here. Let $\E=\O_{\Pone}(1)\oplus
\O_{\Pone}(2)^{\oplus 2}\oplus\O_{\Pone}$, and let $\P=\P(\E)$ be the
$\Pthree$-bundle over $\Pone$. If $t=c_1(\O_{\P}(1))$, then $-K_{\P}\sim
4t-3f$. We can write equations on $\P$ using variables $X,Y,Z,W$ corresponding
to the inclusions of $\O_{\Pone}(1),
\O_{\Pone}(2),\O_{\Pone}(2),\O_{\Pone}$ in $\E$ respectively, and variables
$u,v$ on $\Pone$. Thus
$$T=\{uY^2W^2+vZ^2W^2+X^3W=0\}$$
represents an element of $|-K_P|$, and so $K_T=0$. $T$ is of course reducible,
but is irreducible locally along the line $l=\{X=Y=Z=0\}$. $T$ is also singular
along $l$. Blowing up $l$, it is easy to see that this desingularizes $T$ in a
neighborhood of $l$, and that the exceptional divisor $E$ is a non-normal
rational surface with $E^3=7$. Furthermore, the linear system $|t|$ on $\P$
contracts $l$; this induces a morphism $\pi:T\rightarrow T'$, where $\pi(l)$
is an isolated singular point. This provides an example, again locally, of an
isolated rational Gorenstein point which is the contraction of a non-normal del
Pezzo surface of degree 7.

It is possible to globalize this example by finding a more general element of
$|-K_P|$ with the same local behaviour at $l$, but which is non-singular away
from $l$.

\proclaim Lemma 1.6. Let $\pi:\tilde X\rightarrow X$ be a Type III
contraction with normal exceptional locus, $C\cong\Pone$ the singular locus of
$X$, and suppose $\pi$ is the blow-up of $X$ along $C$. Let
$f:\X\rightarrow\Delta$ be a flat deformation of
$X$ over a one-dimensional disk $\Delta$.
Then either there is a flat deformation $\tilde f:\tilde \X\rightarrow \Delta$
of $\tilde X$, and a map $\pi':\tilde \X\rightarrow \X$ over $\Delta$ which is
a
deformation of $\pi$, or else, after possibly making a finite base-change over
$\Delta$, there is a small projective morphism $\X'\rightarrow\X$ with
$\X'\rightarrow \Delta$ flat such that $\X'_t$ has ${\bf Q}$-factorial terminal
singularities. In particular, if $\pi$ is primitive, in the latter case
$\X'=\X$.

Proof. First suppose that each fibre of $f$ has a
one-dimensional singular locus. Let $S\subseteq\X$ be the reduced two
dimensional
part of the singular locus of $\X$. The central fibre of $S\rightarrow \Delta$
is
the curve
$C\subseteq X$. If we take a general
hyperplane section $\H$ of $\X$, then $\H\rightarrow\Delta$ is a deformation of
a hyperplane section $H$ of $X$, which has a certain number of ordinary double
points. Since as we deform $H$, the number of ordinary double points cannot go
up, we see that $C$ is a reduced fibre of $S\rightarrow\Delta$, and so
$S\rightarrow\Delta$ is a non-singular $\Pone$-bundle over $\Delta$ near
$0\in\Delta$. Thus if we blow up $S$ inside of $\X$, we obtain $\tilde f:
\tilde \X\rightarrow \Delta$, and the proper transform of $X$ is the blow-up of
$X$ along $C$, which is $\tilde X$ by hypothesis. Thus
$\tilde f$ is a deformation of
$\tilde X$, as desired.

Now suppose that the general fibre of $f$ has codimension 3 singularities.
Since $X$ has cDV singularities, so does the general fibre, and so the general
fibre has terminal singularities. Suppose now that the general fibre is not
${\bf Q}$-factorial. Let $D_t$ be a Weil divisor in a general fibre which is
not ${\bf Q}$-Cartier. Possibly after making a base change
$\Delta'\rightarrow\Delta$, we can find a divisor $D$ on $\X$ such that
$D|_{\X_t}=D_t$. $D$ is not ${\bf Q}$-Cartier. Let
$$\R_{\X}(D)=\bigoplus_{m=0}^{\infty} \O_{\X}(mD).$$
If $\R_{\X}(D)$ is a finitely generated $\O_{\X}$-algebra, and we set $\X'
={\bf Proj} (\R_{\X}(D))$, then $\X'$ is normal and $\X'\rightarrow\X$
is
an isomorphism in codimension 1, by [12], Lemma 3.1. Furthermore, the proper
transform of
$D$ is ${\bf Q}$-Cartier on $\X'$. We can repeat this process until $\X'_t$ is
${\bf Q}$-factorial for general $t\in\Delta$. If $\pi$ is primitive,
and $\X_t$ were not ${\bf Q}$-factorial, then since
$\X'_0$ is a partial resolution of $\X_0=X$, $\X'_0=\tilde X$ or $X$. But
$\X'_0\not=\X_0$, so $\X'_0=\tilde X$, which is impossible since then
$\X'\rightarrow
\Delta$ is a flat deformation of $\tilde X$ in which the exceptional divisor
$E\subseteq \tilde X$ does not deform. Thus we obtain a contradiction, and so
$\X'_0=X$ and $\X_t$ is ${\bf Q}$-factorial.

We need then to show that $\R_{\X}(D)$ is finitely generated. We proceed as
follows. Working locally, since $\X$ has only cDV singularities, we can write
$\X$ as a two-dimensional family of du Val surface singularities $\X\rightarrow
S$.
By [12], Lemma 3.2, we can replace $\X$ by any finite cover, and thus, in
particular, by making a base-change, we can assume that the family
$\X\rightarrow
S$ has a simultaneous resolution $\X'\mapright{g}\X\rightarrow S$. By taking a
general fibration
$S\rightarrow S_1$, where $S_1$ is one-dimensional, we get a fibration
$\X'\rightarrow S_1$ of non-singular threefolds. Since threefold flops exist in
families ([15], Theorem 11.10) we can perform flops on $\X'$ until $g^*D$ is
$g$-nef. $g^*D$ then induces a factorization $\X'\rightarrow\X''\mapright{h}
\X$ with $h^*D$ $h$-ample. By
[12], Lemma 3.1, it then follows that
$\R_{\X}(D)$ is finitely generated. $\bullet$

\proclaim Theorem 1.7. Let $\pi:\tilde X\rightarrow X$ be a primitive type
III contraction of a non-singular Calabi-Yau threefold $\tilde X$, contracting
a divisor $E$ to a curve $C\cong \Pone$. Then $X$ is smoothable if $E^3\le 6$.

Proof.
Our goal is to find a deformation $\X\rightarrow\Delta$ of $X$ inducing a map
$\Delta\rightarrow Def(X)$ such that the image of $\Delta$ is not contained
in $\im(Def(\tilde X)\mapright{p} Def(X))$. Here the map $p$, given by
blowing-down deformations, exists by [15], (11.4). If we have such a
deformation, then by Lemma 1.6,
$\X_t$ has
${\bf Q}$-factorial terminal singularities for general
$t\in\Delta$. It then follows from [18] that $\X_t$ is smoothable, so $X$ is
smoothable.

Of course, $Def(\tilde X)$ is smooth by the Bogomolov-Tian-Todorov
unobstructedness theorem. If $Def(X)$ is smooth, then all we need to do is show
that the differential of $p$, $p_*:T^1_{\tilde X}\rightarrow T^1_X$ is a proper
inclusion, where $T^1_{\tilde X}$ and $T^1_X$ are the tangent spaces of
$Def(\tilde X)$ and $Def(X)$ respectively. Since the map $p$ has finite fibres
and we can assume
$\tilde X$ is general in its moduli, we can assume that $p_*$ is injective. We
need to show that $p_*$ is not surjective.
As we will see shortly, however, $Def(X)$ is smooth if $E^3\le 5$, but in the
case $E^3=6$, this may not be the case.

{\it Step 1.} $Def(X)$ is smooth if $E^3\le 5$, and $Def(X)$ is defined
set-theoretically by at most one equation in a neighborhood of the origin of
$T^1_X$ if $E^3=6$.

Proof: Let $\hat X$ be the completion of $X$ along the singular locus
$C\subseteq X$. By [8], Remark 2.7, $Def(X)$ is smooth if
$\ext^2_{\O_X}(\Omega^1_X,\O_{\hat X})=0$. Note that
$\lext^i_{\O_X}(\Omega^1_X,\O_{\hat X})$ is supported on $C$, and
$\lext^2_{\O_X}(\Omega^1_X,\O_{\hat X})=0$ since $X$ has only hypersurface
singularities. Thus, by the local-global spectral sequence for Ext's,
$$\ext^2_{\O_X}(\Omega^1_X,\O_{\hat X})=H^1(\lext^1_{\O_X}(\Omega^1_X,\O_{\hat
X})).$$
Now, locally, if $X\subseteq Y$ with $Y$ four-dimensional and non-singular, we
have the sequence
$$\lhom_{\O_X}(\Omega^1_Y|_X,\F)\rightarrow \lhom_{\O_X}(\I_{X/Y}/\I^2_{X/Y},
\F)\rightarrow \lext^1_{\O_X}(\Omega^1_X,\F)\rightarrow 0$$
for any quasi-coherent sheaf $\F$, or equivalently,
$$\T_Y|_X\otimes \F\rightarrow \N_{X/Y}\otimes \F\rightarrow
\lext^1_{\O_X}(\Omega^1_X,\F) \rightarrow 0.$$
This shows that
$$\lext^1_{\O_X}(\Omega^1_X,\F)\cong \lext^1_{\O_X}(\Omega^1_X,\O_X)\otimes
\F.$$
Thus, since $\lext^1_{\O_X}(\Omega^1_X,\O_X)={\bf T}^1$ is a coherent sheaf
supported on $C$, ${\bf T^1}\cong\lext^1_{\O_X}(\Omega^1_X,\O_{\hat
X})$. By Theorem 1.4 (b), $H^1({\bf T}^1)=0$ if $E^3\le 5$, and
$T^2_{loc}:=\ext^2_{\O_X}(\Omega^1_X,
\O_{\hat X})=H^1({\bf T}^1)$ is one-dimensional if $E^3=6$. The claim in
the case that $E^3=6$ then follows from Theorem 2.2 of [8] and Theorem
1 of [13].
$\bullet$

{\it Step 2.} $\dim\coker(p_*:T^1_{\tilde X}\rightarrow T^1_X)\ge 2$ if $E^3\le
6$.

Proof. By [29], (1.5), the map $p_*$ is the map
$$\ext^1_{\O_{\tilde X}}(\Omega^1_{\tilde X},\O_{\tilde X})\rightarrow
\ext^1_{\O_{\tilde X}}(\pi^*\Omega^1_X, \O_{\tilde X})$$ induced by the natural
map $\pi^*:\Omega^1_X\rightarrow \Omega^1_{\tilde X}$; by applying ${\bf R}
\Gamma$ to both sides of (1.4.3), we see that
$\ext^1(\pi^*\Omega^1_X,\O_{\tilde
X})\cong \ext^1(\Omega^1_X,\O_X)=T^1_X$.

We set $\T'=\lhom_{\O_{\tilde X}}(\pi^*\Omega^1_X,\O_{\tilde X})$.
The morphism
$${\bf R}\Gamma {\bf R}\lhom_{\O_{\tilde X}}(\Omega^1_{\tilde X},\O_{\tilde X})
\rightarrow
{\bf R}\Gamma {\bf R}\lhom_{\O_{\tilde X}}(\pi^*\Omega^1_X,\O_{\tilde
X})$$
induced by $\pi^*\Omega^1_X\rightarrow\Omega^1_{\tilde X}$ induces a morphism
of Grothendieck spectral sequences which gives a diagram
$$
\matrix{H^1(\T_{\tilde X})&\mapright{\cong}&\ext^1_{\O_{\tilde X}}
(\Omega^1_{\tilde X},\O_{\tilde X})&&&&&&\cr
\mapdown{}&&\mapdown{p_*}&&&&&&\cr
H^1(\T')&\hookrightarrow&\hidewidth\ext^1_{\O_{\tilde X}}(\pi^*\Omega^1_X,
\O_{\tilde X})\hidewidth&\mapright{}&&&&\cr
&&&\hidewidth H^0(\lext^1_{\O_{\tilde
X}}(\pi^*\Omega^1_X,
\O_{\tilde
X}))&\mapright{}&
H^2(\T')
&\mapright{\alpha}&\ext^2_{\O_{\tilde X}}
(\pi^*\Omega^1_X,
\O_{\tilde X})\cr}$$
Thus $\dim\coker p_*\ge2 $ if
$$\dim\ker\alpha\le \dim H^0(\lext^1_{\O_{\tilde X}}(\pi^*\Omega^1_X,
\O_{\tilde X})) -2.$$
By Theorem 1.4 c),
$$\dim H^0(\lext^1_{\O_{\tilde X}}(\pi^*\Omega^1_X,
\O_{\tilde X})) \ge 8-E^3+\dim H^0(\O_C(4-E^3)).$$
Thus, since we are assuming that $E^3\le 6$, we just need to show that
$$\dim\ker\alpha\le \dim H^0(\O_C(4-E^3)).$$
We will show this using the sequence (1.4.4). Breaking this up into two short
exact sequences, we have first
$$\exact{\G}{\pi^*\Omega^1_X}{\ker\phi}.\leqno{(1.7.1)}$$
Dualising this shows that $\T'\cong \dual{(\ker \phi)}$, while applying
$\hom_{\O_{\tilde X}}(\cdot,\O_{\tilde X})$ gives us a commutative diagram
$$\matrix{&&\ext^1_{\O_{\tilde X}}(\G,\O_{\tilde X})\cr
&&\mapdown{\delta}\cr
H^2(\dual{(\ker\phi)})&\mapright{\beta}&\ext^2_{\O_{\tilde
X}}(\ker\phi,\O_{\tilde X})\cr
\mapdown{\cong}&&\mapdown{}\cr
H^2(\T')&\mapright{\alpha}&\ext^2_{\O_{\tilde X}}(\pi^*\Omega^1_X,
\O_{\tilde X})\cr}\leqno{(1.7.2)}$$
where $\beta$ also comes from the local-global $\ext$ spectral sequence.
This shows that if $\beta$ is injective, then $\dim\ker\alpha\le
\dim \ext^1_{\O_{\tilde X}}(\G,\O_{\tilde X})$. But by Serre duality,
$\ext^1_{\O_{\tilde X}}(\G,\O_{\tilde X})\cong \dual{H^2(\G)}$. From the proof
of 1.4 c), $\G\cong \pi^*O_C(E^3-4)\otimes\omega_E$, so by Serre duality on
$E$,
$$\eqalign{\dual{H^2(\G)}\cong & H^0(\pi^*\O_C(4-E^3))\cr
\cong & H^0(\O_C(4-E^3)).\cr}$$
Thus we see that if $\beta$ is injective, then $\dim\ker\alpha \le
\dim H^0(\O_C(4-E^3))$ as desired.

To show $\beta$ is injective, we use the other piece of the sequence (1.4.4),
$$\exact{\ker\phi}{\Omega^1_{\tilde X}}{\Omega^1_{\tilde
X/X}}.\leqno{(1.7.3)}$$
Dualising this sequence and using the fact that
$\lext^1_{\O_{\tilde X}}(\Omega^1_{\tilde X/X},\O_{\tilde X})
\cong \lhom_{\O_E}(\Omega^1_{\tilde X/X},\omega_E)$, we obtain a commutative
diagram with exact rows, setting $\T:=\lhom_{\O_E}(\Omega^1_{\tilde
X/X},\omega_E)$,
$$\matrix{ H^1(\T)&\mapright{\gamma_1} &H^2(\T_{\tilde
X})&\mapright{}&H^2(\dual{(\ker\phi)})&\mapright{}&
H^2(\T)\cr
\mapdown{}&&\mapdown{\cong}&&\mapdown{\beta}&&\mapdown{}\cr
\ext^1_{\O_E}(\Omega^1_{\tilde X/X},\omega_E)&\mapright{\gamma_2}&
\ext^2_{\O_{\tilde X}}(\Omega^1_{\tilde X},\O_{\tilde X})&\mapright{}&
\ext^2_{\O_{\tilde X}}(\ker\phi,\O_{\tilde X})&\mapright{}&
\ext^2_{\O_E}(\Omega^1_{\tilde X/X},\omega_E)\cr}\leqno{(1.7.4)}$$
where the vertical maps come from the local-global $\ext$ spectral sequence.
To show that $\beta$ is injective, it will be enough to show the following
claims:

{\it Claim 1:} $H^2(\T)=0$.

{\it Claim 2:} $\im\gamma_1=\im\gamma_2$.

{\it Proof of Claim 1:} Since $\T$
restricted to a general fibre $f$ of $\pi:E\rightarrow C$ is $\O_f$,
$R^1\pi_*\T$ is supported on points.
Thus by the Leray spectral sequence, $H^2(\T)=0$.

{\it Proof of Claim 2:} Using Serre duality, the first square of (1.7.4)
is dual to
$$\matrix{\ext^1_{\O_E}(\dual{(\Omega^1_{\tilde X/X})},\O_E)
&\mapleft{\dual{\gamma_1}}&H^1(\Omega^1_{\tilde X})\cr
\mapup{\delta}&&\mapup{=}\cr
H^1(\Omega^1_{\tilde X/X})&\mapleft{\dual{\gamma_2}}&H^1(\Omega^1_{\tilde
X})\cr
}$$
Clearly $\im\gamma_1\subseteq\im\gamma_2$, and since
$\dim\im\gamma_i=\dim\im\dual{\gamma_i}$, it is enough to show that
$\dim\im\dual{\gamma_2}\le\dim\im\dual{\gamma_1}$. The local-global $\ext$
spectral sequence gives an inclusion
$$H^1(\dual{\dual{(\Omega^1_{\tilde X/X})}})
\rightarrow \ext^1_{\O_E}(\dual{(\Omega^1_{\tilde X/X})},\O_E)$$ through which
$\delta$ factors via the natural map $$\delta':H^1(\Omega^1_{\tilde X/X})
\rightarrow H^1(\dual{\dual{(\Omega^1_{\tilde X/X})}}).$$ Furthermore,
$\pi_*\Omega^1_{\tilde X/X}$ and $\pi_*\dual{\dual{(\Omega^1_{\tilde X/X})}}$
are sheaves of finite length, so
$H^1(\Omega^1_{\tilde X/X})=H^0(R^1\pi_*\Omega^1_{\tilde X/X})$ and
$H^1(\dual{\dual{(\Omega^1_{\tilde
X/X})}})=H^0(R^1\pi_*\dual{\dual{(\Omega^1_{\tilde X/X})}})$. Thus it is enough
to show in the diagram
$$\matrix{H^0(R^1\pi_*\dual{\dual{(\Omega^1_{\tilde
X/X})}})&\mapleft{\gamma_1'}
&H^1(\Omega^1_{\tilde X})\cr
\mapup{\delta'}&&\mapup{=}\cr
H^0(R^1\pi_*\Omega^1_{\tilde X/X})&\mapleft{\dual{\gamma_2}}&
H^1(\Omega^1_{\tilde X})}$$
that $\dim\im \dual{\gamma_2}\le\dim\im\gamma_1'$.
$X$ is ${\bf Q}$-factorial, as follows from [14], Proposition 5-1-6, because
$\pi$ is an extremal divisorial contraction. Thus by
[8], Lemma 4.4
$$\im(H^1(\Omega^1_{\tilde X})\rightarrow H^0(R^1\pi_*\Omega^1_{\tilde X}))
=\im({\bf C}E\rightarrow H^0(R^1\pi_*\Omega^1_{\tilde X})),$$
so $$\im(\dual{\gamma_2})=\im({\bf C}E\rightarrow
H^0(R^1\pi_*\Omega^1_{\tilde X/X})).$$
This image is one dimensional, and the image of $E$ gives a section of
$R^1\pi_*\Omega^1_{\tilde X/X}$ which is clearly supported on all of $C$. Since
$R^1\pi_*\Omega^1_{\tilde X/X}$ and
$R^1\pi_*\dual{\dual{(\Omega^1_{\tilde X/X})}}$ differ on only a finite set of
points, $\delta'(\sigma)$ is non-zero, so $\dim\im\gamma_1'\ge 1$. Thus
$1=\dim\im(\dual{\gamma_2})\le\dim\im\gamma_1'$.

{\it Step 3.} As already observed at the beginning of the proof, Step 2 implies
the desired result if $Def(X)$ is smooth. This is the case if $E^3\le 5$, by
Step 1. If $E^3=6$, then $\dim T^1_X\ge (\dim \im p_*)+2$ by Step 2, and
$Def(X)$ is defined set-theoretically by at most one equation in $T^1_X$. Thus
there still exists a deformation $\X\rightarrow S$ for some $S$ which is at
least one dimensional, such that this deformation does not come from a
deformation of
$\tilde X$, and we finish as before.
$\bullet$

{\it Remark 1.8.} The arguments of Step 1 of the proof of Theorem 1.7 show that
$Def(X)$ can be obstructed if $E^3=6,7$ or $8$. However, if $E^3=7$,
$T^2_{loc}$ is two-dimensional, yet $H^0({\bf T}^1)$ only provides one
additional tangent direction. If $E^3=8$, in fact $H^0({\bf T}^1)=0$. So the
$E^3=6$ case is the only case where one might expect to see an obstructed
deformation space yet still have some actual deformations. Indeed, this is
precisely the case for the original example of a Calabi-Yau threefold with
canonical singularities and obstructed deformations given in [7]. There, the
singular Calabi-Yau threefold has a resolution with exceptional divisor $E$ and
$E^3=6$.

We can get a stronger result if $E^3\le 4$. The following result should be
compared with [8], Theorem 3.8, which implies as a special case that if
$\pi:\tilde
X\rightarrow X$ is a contraction (not necessarily primitive) of a divisor $E$
to a
point and
$E^3\le 4$, then $X$ is smoothable.

\proclaim Theorem 1.9. Let $\tilde X$ be a non-singular Calabi-Yau threefold
and
let $\pi:\tilde X\rightarrow X$ be a contraction of a normal divisor $E$ to a
curve $C\cong\Pone$, and assume that $\pi$ is the blow-up of $X$ along $C$. If
$E^3\le 4$, then there is a partial resolution $X'\rightarrow X$ with $X'$
singular, such that $X'$ is smoothable.

Proof. We first note that the conclusion about $E$ of Theorem 1.1 b) still
holds,
since the proof in [30] only uses the fact that $\tilde X$ is non-singular and
$\pi$ is the blow-up of $X$ along $C$. The proof of Theorem 1.4 b) and c) only
relies on these facts and Theorem 1.1 b). So we can still apply these results.

Now suppose we find a deformation $\X\rightarrow\Delta$ of $X$ inducing a map
$\Delta\rightarrow Def(X)$ such that the image of $\Delta$ is not contained in
$\im(Def(\tilde X)\mapright{p} Def(X))$. Then by Lemma 1.6, we can find, after
making a base-change over $\Delta$, a family $\X'\rightarrow\Delta$ with
$\X'\rightarrow \X$ a small morphism and $\X'_t$ having ${\bf Q}$-factorial
terminal singularities. Set $X':=\X'_0$. By [18], it follows that $\X'_t$,
hence
$X'$, is smoothable. Note that while $X'\rightarrow X$ need not itself be a
small
morphism, $X'$ cannot be non-singular; otherwise the exceptional divisor would
deform in the family $\X'\rightarrow\Delta$, and then $\X'\rightarrow\X$ would
not be a small morphism.

Now, since $E^3\le 4$, Step 1 of the proof of Theorem 1.7 shows that $Def(X)$
is
smooth and so we only need to show that $p_*:T^1_{\tilde X}\rightarrow T^1_X$
is not a surjection. As in Step 2 of the proof of Theorem 1.7, we need to show
that
$$\dim\ker\alpha< \dim H^0(\lext^1_{\O_{\tilde X}}(\pi^*\Omega^1_X,\O_{\tilde
X})).$$
By Theorem 1.4 c),
$$\dim H^0(\lext^1_{\O_{\tilde X}}(\pi^*\Omega^1_X,\O_{\tilde X}))
=length(\F)+ \dim H^0(\O_C(4-E^3)).$$
Since $E^3\le 4$, $\dim H^0(\O_C(4-E^3))>0$, so we just need to show
$$\dim\ker\alpha\le length(\F).$$
Suppose in diagram (1.7.2), the map $\delta$ was zero. Then $\dim\ker\alpha
=\dim\ker\beta$, and it would be enough to show that $\dim\ker\beta\le
length(\F)$.
But by the proof of Theorem 1.4 (c), $\lext^2_{\O_{\tilde
X}}(\ker\phi,\O_{\tilde
X})=0$, and $H^1(\lext^1_{\O_{\tilde X}}(\ker\phi,\O_{\tilde X}))=0$, so the
local-global Ext spectral sequence yields an exact sequence
$$H^0(\lext^1_{\O_{\tilde X}}(\ker\phi,\O_{\tilde X}))\rightarrow
H^2(\dual{(\ker\phi)})\mapright{\beta} \ext^2_{\O_{\tilde
X}}(\ker\phi,\O_{\tilde
X})\rightarrow 0.$$
Again from the proof of Theorem 1.4 c),
$$\dim H^0(\lext^1_{\O_{\tilde X}}(\ker\phi,\O_{\tilde X}))=length(\F),$$
so we see $\dim\ker\beta\le length(\F)$. Thus we just need to show

{\it Claim:} The map $\delta$ in (1.7.2) is zero.

Proof: The Serre dual of $\delta$ is the boundary map
$$H^1(\ker\phi)\mapright{\dual{\delta}} H^2(\G)$$
induced by (1.7.1), so it is enough to show that $\dual{\delta}$ is zero,
or equivalently, that $$H^1(\pi^*\Omega^1_X)\mapright{\xi} H^1(\ker\phi)$$ is
surjective. Now since $H^0(\Omega^1_{\tilde X/X})=H^0(\Omega^1_{E/C})=0$,
(1.7.3)
yields the exact sequence
$$0\rightarrow H^1(\ker\phi) \rightarrow H^1(\Omega^1_{\tilde X}) \rightarrow
H^1(\Omega^1_{E/C}).$$
If we identify $H^1(\Omega^1_{\tilde X})$ with $\Pic \tilde X\otimes {\bf C}$,
this shows that $H^1(\ker\phi)$ is generated by divisors which are all
$\pi$-numerically trivial; i.e. divisors which are pullbacks of divisors on
$X$.
Thus we have a commutative diagram
$$\matrix{H^1(\O_X^*)&&\hidewidth\mapright{\pi^*}\hidewidth
&&H^1(\O^*_{\tilde
X})\cr
\mapdown{dlog}&&&&\mapdown{dlog}\cr
H^1(\Omega^1_X)&\mapright{\pi^*}&H^1(\pi^*\Omega^1_X)&\mapright{}&
H^1(\Omega^1_{\tilde X})\cr}$$
and the image of the composed map $H^1(\O_X^*)\rightarrow H^1(\Omega^1_{\tilde
X})$ generates $H^1(\ker\phi)$. Thus $\xi$ is surjective. $\bullet$

{\hd \S 2. The K\"ahler Cone of Primitive Calabi-Yau Threefolds.}

Having now understood which primitive contractions yield smoothable Calabi-Yau
threefolds, we would like to get a feeling of what kind of constraints this
information imposes on primitive Calabi-Yaus. As a first example, recall the
following theorem of Nikulin:

\proclaim Theorem 2.1. Suppose $X$ is a non-singular Calabi-Yau
threefold. Then either:
\item{(1)} $rank\Pic(X)\ge 41$;
\item{(2)} $X$ has a small contraction;
\item{(3)} There exists a nef ${\bf R}$-divisor $D\in\Pic(X)\otimes {\bf R}$
with
$D^3=0$.

The number 41 is probably far from being optimal. However, the constraints
imposed on primitive Calabi-Yaus easily yield the following optimal
result.

\proclaim Theorem 2.2. Suppose $X$ is a non-singular primitive Calabi-Yau
threefold. Then either:
\item{(1)} $rank\Pic(X)=1$;
\item{(2)} $X$ has a small contraction;
\item{(3)} There exists a nef ${\bf R}$-divisor $D\in\Pic(X)\otimes {\bf R}$
with
$D^3=0$.

Proof: Suppose that neither case (2) or (3) occur. By [30], the closure of
the K\"ahler cone $\bar\K$ of $X$ is rational polyhedral, and every codimension
one face induces a primitive divisorial contraction. Following the notation of
[19], let $R_1,\ldots, R_m$ be the extremal rays of $NE(X)$, the cone of
effective curves, and denote by
$D(R_i)$ the divisor contracted by contracting the extremal ray $R_i$. By
[19], Theorem 1.3.2 (1), the $D(R_i)$ are all distinct divisors.
Now let $\{R_i\}_{i\in I}$ be an $E$-set. (See [19], \S 1.1 for the definition
of
an $E$-set.) Then by Theorem 1.2.13 of [19], $\{R_i\}_{i\in I}$ satisfies
condition (iii) of \S 1.1 of [19], so there is a non-zero nef divisor
$D=\sum_{i\in I} a_iD(R_i)$ with $a_i\ge 0$. Now according to Theorem 0.4 and
Theorem 5.8 of [8], $D(R_i)$ is rational and satisfies $D(R_i)^3\ge 7$;
otherwise
the contraction induced by $R_i$ would yield a smoothable singular Calabi-Yau
threefold. It then follows as in the proof of [8], Theorem 5.2 that
$c_2(X).D(R_i)=12-2D(R_i)^3<0$. But then $c_2(X).D<0$, but this
is impossible, since $c_2$ is non-negative on the K\"ahler cone of $X$.
$\bullet$

Let's speculate about each case in turn.

(1) Every Picard number 1 Calabi-Yau is obviously primitive. Some completely
new idea will be needed to understand such threefolds. I hope that these
threefolds will prove to form a more manageable class.

(2) There are examples of primitive Calabi-Yau threefolds with small
contractions.

{\it Example 2.3.} Take a general anti-canonical hypersurface in the
$\Pthree$-bundle $\P(\O_{\Pone}(-1)\oplus\O_{\Pone}^{\oplus 2}\oplus
\O_{\Pone}(1))$ over $\Pone$. (See [8], [27] for details on this threefold.)
Such a Calabi-Yau has Picard number 2, and the two boundaries of the K\"ahler
cone correspond to a K3 fibration and a small contraction of a single $\Pone$
to
an ordinary double point. However, if we flop this $\Pone$, we obtain a
non-primitive Calabi-Yau which has a contraction of a rational ruled surface
$E$ to a $\Pone$ with $E^3=-2$.

So we might ask for yet a stronger condition:

\proclaim Definition 2.4. A Calabi-Yau threefold $X$ is {\it birationally
primitive} if every minimal model of $X$ is primitive.

I do not know the answer to

\proclaim Question 2.5. Can a birationally primitive Calabi-Yau threefold have
more
than one minimal model?

I would hesitate to turn this into a conjecture at this point; it could merely
reflect an ignorance of examples.

(3) is the case we know the most about. It appears to be quite a reasonable
conjecture that if there exists a nef ${\bf R}$-divisor $D$ with $D^3=0$, then
there is a nef ${\bf Q}$-divisor with the same property. This would follow from
some cone conjectures: see [33] for a survey of these conjectures. If
$D$ is a nef divisor with $D^3=0$, and $D.c_2>0$, then $|nD|$ induces a K3 or
elliptic fibration for sufficiently large $n$ (see [20]). If
$D.c_2=0$,  then we don't know yet whether $|nD|$ will always induce a
fibration (elliptic or abelian), but see [32] for results in this direction.
In general one should expect to obtain such a fibration.

Once we know a Calabi-Yau has a fibration, this gives us a great deal of
information. For example, up to birational equivalence, there are only a
finite number of elliptic Calabi-Yau threefolds [6]. Given the classification
methods known for elliptic fibrations, I believe it would be a quite tractable
task to find all primitive elliptic Calabi-Yau fibrations. Perhaps K3 and
abelian surface fibrations will prove to be similarly tractable.

{\it Example 2.6.} It is often harder to tell if a Calabi-Yau is primitive than
to show one is not primitive. Some Calabi-Yau threefolds are obviously
primitive,
as mentioned in the introduction: if a Calabi-Yau has no birational
contractions,
then it is primitive. This includes examples such as hypersurfaces of bidegree
$(3,3)$ in $\Ptwo\times\Ptwo$, double covers of $\Pone\times\Ptwo$ and double
covers of $\Pone\times\Pone\times\Pone$, as well as any Picard number 1
threefold. Example 2.3 yields a primitive Calabi-Yau which does have a
birational
contraction. However, if the Picard number is much larger than 2, then there
could be many birational contractions, most of them not primitive contractions,
and at the moment we can't really check each of these models for smoothability.
However, we can give some likely candidates for primitive Calabi-Yaus with
larger
Picard number.

Let $S$ be one of the rational scrolls $F_{12}$, $F_8$, $F_6$ or $F_4$. Let
$X_{12},X_8,X_6$ and $X_4$ be elliptic fibrations over each of these surfaces
respectively given by Weierstrass equation
$$Y^2=X^3+aX+b$$ with general
$a\in\Gamma(S,\omega_S^{\otimes -4}), b\in\Gamma(S,\omega_S^{\otimes -6})$.
It is easy to check that these are Calabi-Yau threefolds, and have non-singular
minimal models. Let $C_0$ be the negative section of $S$. Then over $C_0$, the
fibration $X_{12}\rightarrow S$ has fibres of Kodaira type $II^*$, the
fibration
$X_8\rightarrow S$ has fibres of Kodaira type $III^*$, the fibration
$X_6\rightarrow S$ has fibres of Kodaira type $IV^*$, and the fibration
$X_4\rightarrow S$ has fibres of Kodaira type $I_0^*$, and in each case the
inverse image of $C_0$ is a union of minimal rational scrolls. We have
$\rho(X_{12})=11$, $\rho(X_{8})=10$, $\rho(X_6)=9$ and $\rho(X_4)=7$.
Furthermore, each of these Calabi-Yaus has only one minimal model. Every
primitive contraction is a type III contraction of a minimal rational scroll,
which will be either one of the scrolls in the inverse image of $C_0$ or the
unique section of the fibration. It seems very likely that each of
these Calabi-Yaus is primitive. $X_{12}$ is of special interest since
$\chi(X_{12})=-960$, and this is apparently the smallest known Euler
characteristic of a Calabi-Yau threefold. If $X_{12}$ were not primitive, we
would
expect to find a Calabi-Yau with even smaller Euler characteristic.

{\hd \S 3. Relating primitive Calabi-Yau threefolds to general Calabi-Yau
threefolds.}

Our goal is to shed some light on Questions 0.3, and their relationship to the
classification of primitive Calabi-Yaus. In particular, we would like to
explore
possible corollaries to the conjecture that there are only a finite number of
families of primitive or birationally primitive Calabi-Yaus. We first need to
consider the question of deformation invariance of the K\"ahler cone for
Calabi-Yau
threefolds with canonical singularities. In particular, let
$\X\rightarrow\Delta$ be a one-parameter smoothing of a Calabi-Yau $X$ with
canonical singularities. We would like to compare the K\"ahler cone of $\X_0$,
$\bar\K_0$, with the K\"ahler cone of $\X_t$, $\bar \K_t$, for general $t$.
Now, in general, the Picard numbers $\rho(\X_0)$ might be smaller than
$\rho(\X_t)$. For example, if
$X=\X_0$ is obtained by contracting an elliptic scroll on $\tilde X$ to a
curve, then $\rho(\X_t)=\rho(\tilde X)$, while $\rho(X)=\rho(\tilde X)-1$.
The problem here is that we contracted ``too much''. We can deform $\tilde X$
to the same threefolds we can deform $X$ to. So there was no need to contract
the elliptic scroll. More generally, if we wish to make a birational
contraction $\pi:\tilde X\rightarrow X$ and then smooth $X$, then we should
choose $\pi$ to be a ``minimal'' such contraction which yields the desired
smoothing. We make this precise below.

\proclaim Proposition 3.1. Let $f:\X\rightarrow\Delta$ be a
deformation of a Calabi-Yau $X$ with canonical singularities, and suppose that
$f:\X\rightarrow\Delta$ satisfies the condition
$$\vbox{For $p:S\rightarrow\Delta$ any morphism, $\X_S=\X\times_{\Delta} S$,
and for any Weil divisor $D$ on $\X_S$ which is Cartier outside finitely many
fibres, $D$ is ${\bf Q}$-Cartier.}\leqno{(3.1.1)}$$
Then, shrinking
$\Delta$ if necessary, the Picard number of $\X_t$ is constant for
$t\in\Delta$,
and
$(R^1f_*\O_{\X}^*)\otimes_{\boldz} {\bf Q}\cong R^2f_*{\bf Q}$ is a constant
sheaf on
$\Delta$.

Proof. If $f:\X\rightarrow\Delta$ satisfies the hypotheses, then it satisfies
conditions (12.2.1.1)-(12.2.1.3) of [15]. The result then follows from [15],
12.2.5 and the fact that $H^1(\O_{\X_t})=H^2(\O_{\X_t})=0$ for all
$t\in\Delta$. We have to tensor with ${\bf Q}$ as the specialization of a
Cartier divisor might only be ${\bf Q}$-Cartier.
$\bullet$

Given an arbitrary smoothing $f:\X\rightarrow\Delta$, can we modify it so that
it satisfies condition (3.1.1)? For this, we need

\proclaim Definition. Let $X$ be a normal variety. A ${\bf Q}$-factorialization
$\pi:\bar X\rightarrow X$ is an isomorphism in codimension one with $\bar X$
normal such that $\bar X$ is ${\bf Q}$-factorial.

\proclaim Proposition 3.2. Let $f:\X\rightarrow\Delta$ be a
one-parameter smoothing of a Calabi-Yau $X$. Then there is a finite morphism
$\Delta'\rightarrow\Delta$ such that if $\bar\X\rightarrow\X_{\Delta'}$ is a
${\bf Q}$-factorialization then $\bar\X\rightarrow \Delta'$ satisfies condition
(3.1.1).

Proof. Shrinking $\Delta$ if necessary, $R^1f_*\O_{\X}^*|_{\Delta-0}$ is a
local system which by [15], 12.2.5 has finite monodromy. Thus, we can make a
finite base-change $\Delta'\rightarrow\Delta$ ramified only at $0$ such that
$R^1f_*\O_{\X_{\Delta'}}^*|_{\Delta'-0}$ is a constant local system.
It is then clear that if $\bar\X\rightarrow\X_{\Delta'}$ is a ${\bf
Q}$-factorialization of $\X_{\Delta'}$, then $\bar\X\rightarrow\Delta$
satisfies (3.1.1). $\bullet$

If $f:\X\rightarrow\Delta$ is a family satisfying (3.1.1), then we can compare
the K\"ahler cones in the family. If $f$ is a smoothing of $X=\X_0$, then by
[30] we can assume that $\bar\K_t$ is constant for $t\not=0$, by shrinking
$\Delta$ if necessary. We then have $\bar\K_0\subseteq\bar \K_t$. Now
we know from [30] and [17] that equality can fail to hold if $\X_0$
contains certain types of ruled surfaces. However, when $\X_0$ has canonical
singularities, the situation is worse: even small contractions in $\X_0$ may
not deform to the general fibre.

{\it Example 3.3.} Let $P_1=\P(\O_{\Pone}^{\oplus 4})=\Pone\times\Pthree$ and
let $P_2=\P(\O_{\Pone}(-1)\oplus\O_{\Pone}^{\oplus 2}\oplus\O_{\Pone}(1))$. Let
$C\subseteq P_2$ be the section of $P_2$ corresponding to the surjection
$\O_{\Pone}(-1)\oplus\O_{\Pone}^{\oplus 2}\oplus\O_{\Pone}(1)\rightarrow
\O_{\Pone}(-1)$. By [7], it is possibly to construct a family $\X\rightarrow
\Delta$ where $\X_0$ is an anticanonical hypersurface in $P_2$ with double
points along $C$, and $\X_t$ for $t\not=0$ is a non-singular anticanonical
hypersurface in $P_1$. If we write $\Pic\X_t=\boldz T\oplus \boldz F$, where
$T$ is the restriction of $c_1(\O_{P_i}(1))$ to $\X_t$ and $F$ the restriction
of a fibre of $P_i\rightarrow\Pone$, then the K\"ahler cone of $\X_0$ is
spanned by $F$ and $T+F$, the latter divisor contracting $C$, while the
K\"ahler cone of $\X_t$, $t\not=0$ is spanned by $T$ and $F$. This shows that
the K\"ahler cone might fail to be constant because a small contraction fails
to deform.

To rectify this problem, we need the existence of fourfold flops. Recall

\proclaim Definition. Let $X^-$ be a normal variety and $D$ a ${\bf Q}$-Cartier
divisor on $X^-$, and let $\phi:X^-\rightarrow Y$ be a birational contraction
which is an isomorphism in codimension one with $-D$ $\phi$-ample and $K_X$
$\phi$-trivial. The $D$-flop of $\phi$ is a normal variety $X^+$ and a diagram
$$\matrix{X^-&&\cdot\cdot\mapright{\psi}&&X'\cr
&\mapse{\phi}&&\mapsw{\phi^+}&\cr
&&Y&&\cr}$$
with $\psi$ a birational map, $\phi^+$ an isomorphism in codimension one and
$\psi(D)$
$\phi^+$-ample.

We note that existence of flops for fourfolds is not known. However, it is
known ([16]) that there is no infinite sequence of $D$-flops on a fourfold
with terminal singularities.

If fourfold flops exist, however, then we can proceed as follows. If $D$
is a divisor on $\X$, with $D_t\in\bar\K_t$, $t\not=0$, then there is a family
$\X'\rightarrow \Delta$ related by flops over $\Delta$ to $\X$ with $\X_0'$
birationally equivalent to $\X_0$, such that $D$ is nef. We use this for

\proclaim Theorem 3.4. Suppose ${\bf Q}$-factorializations and flops exist for
fourfolds with terminal singularities. Then every Calabi-Yau
threefold
$\tilde X$ with ${\bf Q}$-factorial terminal singularities has a minimal model
which is a crepant ${\bf Q}$-factorial terminal resolution of a Calabi-Yau
$X$, with
$X$ the central fibre of a family $\X\rightarrow\Delta$ and $\X_t$ a
non-singular birationally primitive Calabi-Yau for $t\not=0$.

Proof. We proceed by induction on $\rho(\tilde X)$. If $\tilde X$ is
birationally primitive, we are done. This is in particular the case if
$\rho(\tilde
X)=1$.

Suppose $\tilde X$ is a non-singular, non-birationally primitive Calabi-Yau.
Then, changing the minimal model of $\tilde X$ if necessary, there is a
contraction
$\tilde X\rightarrow X$ and a smoothing $\X\rightarrow\Delta$ of $X$. By
Proposition 3.2, after making a base-change, we can assume
$\X\rightarrow\Delta$ satisfies (3.1.1), and we can then replace $\X$ with its
${\bf Q}$-factorialization. This might change $\X_0$, but it will still be
birational to $\tilde X$. We replace $X$ by $\X_0$. By Proposition 3.1,
$\rho(\X_t)=\rho(\X_0) =\rho(X)<\rho(\tilde X)$ for general $t\in\Delta$. If
$\X_t$ is birationally primitive, we are again done. If not, by induction there
is a minimal model $\X_t'$ of
$\X_t$, a family $\Y\rightarrow\Delta'$ with $\X_t'\rightarrow \Y_0$ a crepant
resolution, with $\Y_t$ birationally primitive. Let $H$ be an ample divisor on
$\Y$, which we pull back to $H\in\Pic\X_t'$, where $H$ is nef. Identifying
$\Pic\X_t$ with
$\Pic
\X_t'$, we denote also by
$H\in\Pic\X_t$ the proper transform of $H$ on $\X_t$. $H$ is in
the moving cone of $\X_t$ (see [12] for the definition of the moving cone). By
Proposition 3.1,
$H$ defines a ${\bf Q}$-Cartier Weil divisor on
$\X$. Since $\X_0$ has canonical singularities and $\X_t$ is smooth, $\X$ has
terminal singularities by [28]. By performing fourfold flops on $\X$ over
$\Delta$,  we can find a family
$\X'\rightarrow
\Delta$ birationally equivalent to $\X\rightarrow\Delta$ such that $H$ is nef
on $\X'$. $H$ then induces a contraction $\X'\rightarrow \X''$ over $\Delta$
such that
$\X''_t\cong\Y_0$. Since $\Y_0$ can be deformed to a birationally primitive
Calabi-Yau,
$\X_0''$ can also be deformed to a birationally primtive Calabi-Yau. This gives
the desired family.

If $\tilde X$ is singular with ${\bf Q}$-factorial terminal singularities, then
by
[18], we can smooth
$\tilde X$, obtaining a family $\X\rightarrow\Delta$ whose general fibre is
smooth
and $\X_0=\tilde X$. We then proceed as before. $\bullet$

This answers Question (1) of the introduction, assuming the existence of a
finite
classification of (birationally) primitive Calabi-Yaus and an understanding of
fourfold birational geometry. What are the obstacles to answering the other
questions given this information?

First, in passing from Question (1) to Question (2), let us try to apply the
same
proof as in Theorem 3.4. We immediately run into trouble, since when we pass to
the ${\bf Q}$-factorialization, $\tilde X$ may no longer be a resolution of
$\X_0$. Secondly, when we perform flops on $\X$, this might also change the
birational model of $\X_0$, and it may no longer have $\tilde X$ as a
resolution.
Thus some better understanding of the birational geometry involved is
necessary.

We now discuss Question (3).
Fix a non-singular Calabi-Yau threefold $X$. We would like to know that $X$ has
a
finite number of degenerations. In other words, we would like a scheme $\SS$ of
finite type, and a flat family $\X\rightarrow\SS$ of Calabi-Yau threefolds with
canonical singularities such that for any family $\X'\rightarrow\Delta$ with
$\X'_t$ deformation equivalent to $X$, $\X'_0$ is birational to $\X_s$ for some
$s\in\SS$. To construct $\SS$, choose an ample divisor $H$ on $X$. Suppose that
there is an integer $r$ satisfying the condition
$$\vbox{For every family $f:\X'\rightarrow\Delta$ with $\X'_t$ deformation
equivalent to $X$ and $\X'_0$ having canonical singularities and such that the
divisor $H$ on $\X'_t$ extends to a ${\bf Q}$-Cartier $f$-ample divisor $\H$ on
$\X'$, $r\H$ is Cartier.}\leqno{(3.5)}$$
If this is the case, then we can just take $\SS$ to be a suitable open subset
of
the irreducible component of the Hilbert scheme containing a point
corresponding to
$X$ embedded via
$10rH$, and $\X\rightarrow\SS$ the universal family. Then given any family
$f:\X'\rightarrow\Delta$, with $\X'_t$ deformation equivalent to $X$, we can,
after making a base-change, a ${\bf Q}$-factorialization and some flops, assume
that $\X'$ has an $f$-ample ${\bf Q}$-Cartier divisor $\H$ extending $H$. Then
$r\H$ is Cartier, and by [21], $10r\H$ is $f$-very ample, and then $\X_0'=\X_s$
for
some $s\in\SS$.

Unfortunately, it is not clear when (3.5) holds. $\X$ has fourfold terminal
Gorenstein singularities. [22] gives a sequence of fourfold terminal
Gorenstein singularities with
${\bf Q}$-Cartier Weil divisors $D$ such that the minimum $r$ such that $rD$ is
Cartier is unbounded. So it may be possible that for larger and larger values
of
$r$ in (3.5), one will get more and more possible degenerations. Thus, in
particular, even if there are only a finite number of families of primitive
Calabi-Yau threefolds, it still may be possible that they might have an
infinite
number of degenerations, causing problems for Questions (3) and (4).

One possible way around this problem could come from a greater understanding of
the path one follows from an arbitrary Calabi-Yau to a primitive Calabi-Yau.
For
example, in \S 2, we obtained the combinatorial restrictions of Theorem 2.2
only
by analysing smoothability of primitive contractions. Suppose every
non-primitive
Calabi-Yau threefold $\tilde X$ has a primitive contraction $\pi:\tilde
X\rightarrow X$ with $X$ smoothed by a family $\X\rightarrow\Delta$. It is easy
to understand the index of a ${\bf Q}$-Cartier divisor on $\X$ since $\X$ will
have relatively simple terminal singularities, and so we would be able to get
around the problem of (3.5) by limiting attention to such simpler
degenerations.
But at this point this is mere speculation. A deeper analysis of the K\"ahler
cone
of primitive Calabi-Yau threefolds or a deeper analysis of degenerations of
Calabi-Yau threefolds is needed.

{\hd Bibliography}

\item{[1]} Candelas, P., Green, P., and H\"ubsch, T., ``Rolling Among
Calabi-Yau
Vacua,'' {\it Nucl. Phys.} {\bf B330} (1990), 49--102.
\item{[2]} Friedman, R., ``Simultaneous Resolution of Threefold Double
Points,'' {\it Math. Ann.,} {\bf 274}, (1986) 671--689.
\item{[3]} Green, P., and H\"ubsch, T., ``Possible Phase Transitions Among
Calabi-Yau Compactifications,'' {\it Phys. Rev. Lett.} {\bf 61}, (1988),
1163--1166.
\item{[4]} Green, P., and H\"ubsch, T., ``Connecting Moduli Spaces of
Calabi-Yau
Threefolds,'' {\it Comm. Math. Phys.} {\bf 119}, (1988), 431--441.
\item{[5]} Greene, B., Morrison, D., and Strominger, A., ``Black Hole
Condensation and the Unification of String Vacua,'' {\it Nucl. Phys.} {\bf
B451},
(1995) 109-120.
\item{[6]} Gross, M., ``A Finiteness Theorem for Elliptic Calabi-Yau
Threefolds,'' {\it Duke. Math. J.,} {\bf 74}, (1994) 271--299.
\item{[7]} Gross, M., ``The Deformation Space of Calabi-Yau $n$-folds with
Canonical Singularities can be Obstructed,'' To appear in {\it Essays in Mirror
Symmetry II}.
\item{[8]} Gross, M., ``Deforming Calabi-Yau Threefolds,'' preprint (1995).
\item{[9]} Hartshorne, R., {\it Residues and Duality,} LNM 20, Springer-Verlag
1966.
\item{[10]} Hartshorne, R., {\it Algebraic Geometry,} Springer-Verlag 1977.
\item{[11]} Kawamata, Y., ``Minimal Models and the Kodaira Dimension
of Algebraic Fiber Spaces,'' {\it J. Reine Angew. Math,} {\bf 363}, (1985)
1--46.
\item{[12]} Kawamata, Y., ``The Crepant Blowing-up of 3-dimensional Canonical
Singularities and its Application to the Degeneration of Surfaces,'' {\it Ann.
of
Math.,} {\bf 127}, (1988) 93--163.
\item{[13]} Kawamata, Y., ``Unobstructed Deformations II,'' {\it J. Algebraic
Geometry,} {\bf 3}, (1995) 277--280.
\item{[14]} Kawamata, Y., Matsuda, K., and Matsuki, K., ``Introduction to the
Minimal Model Program,'' {\it Algebraic Geometry, Sendai, Adv. Stud. Pure
Math.}
Vol. 10, T. Oda, ed., Kinokuniya-North-Holland, 1987, 283--360.
\item{[15]} Koll\'ar, J., and Mori, S., ``Classification of Three-dimensional
Flips,'' {\it J. Amer. Math. Soc.}, {\bf 5}, (1992) 533--703.
\item{[16]} Matsuki, K., ``Termination of Flops for 4-folds,'' {\it Amer. J.
Math.}, {\bf 113}, (1991) 835--859.
\item{[17]} Namikawa, Y., ``On Deformations of Calabi-Yau 3-folds with
Terminal
Singularities,'' {\it Topology}, {\bf 33} (1994), 429--446.
\item{[18]} Namikawa, Y. and Steenbrink, J., ``Global Smoothing of Calabi-Yau
three-folds,'' to appear in {\it Inv. Math}.
\item{[19]} Nikulin, V., ``Diagram Method for 3-folds and its Application to
K\"ahler Cone and Picard Number of Calabi-Yau 3-folds I,'' preprint (1993).
\item{[20]} Oguiso, K., ``On Algebraic Fiber Space Structures on a Calabi-Yau
3-fold,'' {\it Int. J. Math.,} {\bf 4}, (1993) 439--465.
\item{[21]} Oguiso, K., and Peternell, T., ``On Polarized Canonical Calabi-Yau
Threefolds,'' {\it Math. Ann.}, {\bf 301}, (1995) 237--248.
\item{[22]} Onn, S., and Sturmfels, B., ``A Note on Lattice Simplices and Toric
Varieties,'' {\it Amer. J. Math.,} {\bf 116}, (1994) 1337--1339.
\item{[23]} Ran, Z., ``Deformations of Maps,'' in {\it Algebraic Curves and
Projective Geometry,} E. Ballico, C. Ciliberto, Eds. LNM 1389, Springer-Verlag
(1989).
\item{[24]} Reid, M., ``Canonical 3-folds,'' {\it G\'eom\'etrie Alg\'ebrique
Angers,} (A. Beauville, ed.), Sijthoff+Noordhoof, (1980) 273--310.
\item{[25]} ``The Moduli Space of 3-folds with $K=0$ May Nevertheless Be
Irreducible,'' {\it Math. Ann.,} {\bf 278}, (1987) 329--334.
\item{[26]} Rotman, J., {\it An Introduction to Homological Algebra,} Academic
Press, (1979).
\item{[27]} Ruan, Y., ``Topological Sigma Model and Donaldson Type Invariants
in
Gromov Theory,'' preprint (1993).
\item{[28]} Stevens, J., ``On Canonical Singularities as Total Spaces of
Deformations,'' {\it Abh. Math. Sem. Univ. Hamburg}, {\bf 58}, 275--283,
(1988).
\item{[29]} Wahl, J., ``Equisingular Deformations of Normal Surface
Singularities, I,'' {\it Annals of Math.,} {\bf 104} (1976), 325--356.
\item{[30]} Wilson, P.M.H., ``The K\"ahler Cone on Calabi-Yau Threefolds,''
{\it Invent. Math.,} {\bf 107}, (1992) 561--583.
\item{[31]} Wilson, P.M.H., {\it Erratum} to ``The K\"ahler Cone on Calabi-Yau
Threefolds,''
{\it Invent. Math.,} {\bf 114}, (1993) 231--233.
\item{[32]} Wilson, P.M.H., ``The Existence of Elliptic Fibre Space Structures
on
Calabi-Yau Threefolds,'' {\it Math. Ann.}, {\bf 300}, (1994) 693--703.
\item{[33]} Wilson, P.M.H.,  ``The Role of $c_2$ in Calabi-Yau Classification:
A
Preliminary Survey,'' To appear in {\it Essays in Mirror Symmetry II}.
\item{[34]} Wilson, P.M.H., ``Symplectic Classes on Calabi-Yau Threefolds,''
preprint, 1995.

\end